\input amstex
\documentstyle{amsppt}

\document

\topmatter
\title
Complex horospherical transform on real sphere
\endtitle

\author Simon Gindikin \endauthor
\dedicatory To Francois Treves  \enddedicatory

\address Departm. of Math., Hill Center, Rutgers University,
110 Frelinghysen Road, Piscataway, NJ 08854-8019
\endaddress
\email gindikin\@math.rutgers.edu \endemail \abstract We define a
new integral transform on the real sphere which is invariant
relative to the orthogonal group and similar to the horospherical
Radon transform for the hyperbolic space. This transform involves
complex geometry associated with the sphere.
\endabstract

\endtopmatter

In integral geometry on hyperbolic spaces and other noncompact
symmetric spaces there are 2 versions of the Radon transform:
geodesic and horospheric \cite {GGG03}. The horospherical Radon
transform has more essential connections with the harmonic
analysis on these spaces than its geodesic analogue. The geodesic
version of the Radon transform is well known on the sphere. It is
the famous Minkowsky-Funk transform of integration along the
big subspheres \cite {GGG03}. Let us recall that this transform
was discovered earlier than the Radon transform. We want in this
note to construct the analogue of the horospherical transform on
the sphere. At first glance it looks strange since there are no
horospheres on the sphere, but we will show that such a transform
exists if in the geometrical background we replace real
horospheres by complex ones. We will follow the idea of \cite
{Gi00} which was developed for some pseudohyperbolic spaces, but, as
we will see, gives an interesting construction on the sphere
also. Usually, in integral geometry we integrate functions (or
other analytic objects) along some submanifolds. In other words we
integrate along a manifold $\delta$-functions with supports on
these manifolds. The idea is to replace these $\delta$-functions with
Cauchy kernels with singularities on similar complex submanifolds,
and then instead of real submanifolds, we are interested in
corresponding complex submanifolds which have no real points. Such
complex horospheres exist for the real sphere and we can
define the horospherical transform. Such applications of the
complex geometry to analysis on real manifolds can give significant
advantages; they are very much in the spirit of the classical
geometry of XIX century.

Let $S=S^{n-1}$ be the $(n-1)$-dimensional sphere in $\Bbb R^n$:

$$\Delta(x)=x^2_1+\cdots +x_n^2=1$$
and $\Bbb C S=\Bbb CS^{n-1}$ be its complexification in $\Bbb C^n:
\Delta(z)=1, z=x+iy $. Let $\Xi$ be the cone
$\Delta(\zeta)=0,\zeta=\xi+i\eta\neq 0$. Let us consider $\Bbb
C^n_z$ and $\Bbb C^n_\zeta$ as dual relative to the form

$$ \zeta \cdot z=\zeta_1z_1+\cdots +\zeta_nz_n.$$

We will call complex horospheres sections $\Cal H(\zeta)$ of the
complex sphere $\Bbb C S$ by the isotropic hyperplane:

$$\Cal H(\zeta)=\{z\in \Bbb C S; \zeta \cdot z=1\}, \quad \zeta \in \Xi.$$
We will also sometimes parameterize horospheres by homogeneous
coordinates $(\zeta,p)$:

$$\Cal H(\zeta,p)=\{z\in \Bbb CS;\zeta \cdot z=p\}, \quad \zeta \in \Xi,p \in \Bbb C,p\neq 0,$$
where $\Cal H(\lambda \zeta,\lambda p)=\Cal H (\zeta
,p),\lambda\neq 0.$ We do not consider the degenerate horospheres
with $p=0$. The cone $\Xi$ in real coordinates is defined by the
conditions
$$\Delta (\xi)=\Delta (\eta),\qquad \xi \cdot \eta=0.$$
\proclaim {Lemma} The horosphere $\Cal H(\zeta)$ does not
intersect the real sphere $S$ if and only if
$$0<\Delta(\xi)=\Delta (\eta)<1.$$\endproclaim It is a simple direct computation in which it is convenient to use
the rotation invariancy. Let us denote the domain in the cone
$\Xi$ defined in Lemma, through $\Xi_+$. In homogeneous
coordinates this condition is $0<\Delta (\xi)<|p|^2$. $\Xi_+$ is a
Stein submanifold in $\Xi$.

The boundary $\partial \Xi_+$ admits a natural fibering over $S$
on $(n-2)$-dimensional spheres $$\gamma_1(x)=\{\zeta=x+i\eta, \Delta
(\eta)=1,x\cdot \eta=0\},x\in S.$$ The corresponding horospheres
intersect $S$ in one point $x$.

In this note we will assume that the functions $f\in C^\infty(S)$. Of
course, it is possible to define the horospherical transform under
very weak conditions and it would be interesting to consider
different  Paley-Wiener theorems for this transform. We will need
some notations for differential forms. Let us denote through
$[a_1,\dots,a_n]$ the determinant of the matrix with the columns
$a_1,\dots,a_n$ some of which can be 1-forms. We expand such
determinants from left to right and use the exterior product for
the multiplication of 1-forms. Such a determinant with identical
columns can differ from zero: $[dx,\dots,dx]=n! dx_1\wedge\cdots
\wedge dx_n$.We will write $a^{\{k\}}$ if a column $a$ repeats $k$
times.

We will use the interior product of forms $\varphi \rfloor\psi$.
It is a form $\alpha$ such that $ \varphi \wedge\alpha=\psi$; its
restriction on a submanifold where $\varphi =0$ is uniquely
defined. If $\varphi=df$ where $f$ is a function then $df\rfloor
\psi$ up to a multiplicative constant is the residue of $\psi/f$
on $\{f=0\}$.

 Let us define the
horospherical transform as

$$\hat f(\zeta,p)=\int _S \frac {f(x)}{\zeta\cdot x -p} [x,dx^{\{n-1\}}],\quad \zeta \in \Xi_+,$$
where $[x,dx^{\{n-1\}}]=2d(\Delta(x))\rfloor
[dx^{\{n\}}]=(n-1)!\sum _{1\leq j\leq n} (-1)^{(j-1)}x_j
\bigwedge_{i\neq j} dx_i $ is the invariant measure on $S$.

 We have
$$\hat f(\lambda \zeta .\lambda p)=\lambda^{-1} \hat f(\zeta, p).$$  Let
$\hat f(\zeta)=\hat f(\zeta, 1)$. Apparently, the transform is
well defined and $\hat f(\zeta)$ is the holomorphic function in
$\Xi_+$. Also under our conditions on $f$ boundary values of $\hat
f$ are well defined. Our aim is to find an inversion formula for
the horospherical transform $f\rightarrow \hat f$.

\proclaim {Theorem} There is an inversion formula

$$\gather f(x)=\frac{n-1}{2(-2\pi)^{n-1}}\int_{\gamma_1(x)} \Cal L_p\hat f(x+i\eta,p)|_{p=1} [x,\eta, d\eta
^{\{n-2\}}],\\\Cal L_p=\frac {n-2}p \frac {\partial^{(n-3)}
}{\partial p^{(n-3)}} -2 \frac {\partial^{(n-2)} }{\partial
p^{(n-2)}} .\endgather
$$\endproclaim
Here we take the boundary values of $\hat f$ on $\partial \Xi$ and
integrate along cycles $\gamma_1(x).$

The proof of this theorem follows the ideas of \cite {Gi95,Gi
00} (see also Ch.5 in \cite {GGG03}). We will construct the
decomposition on plane waves on the sphere. Namely, let us
consider the differential form

$$\kappa_x[f]=\frac {f(u)}{(\zeta \cdot (u-x))^{n-1}}
[u,du^{\{(n-1\}}]\wedge [u+x,\zeta,d\zeta ^{\{n-2\}}],\qquad u\in
S,\zeta \in \Bbb C^n\backslash 0.$$ Here $x\in S$ is a fixed
point. It is the analogue for the sphere $S$ of the form in
decomposition on plane waves in $\Bbb R^n$. It is important that
we permit complex $\zeta$. The crucial point is that this form is
closed. It is the direct consequence of a technical lemma \cite
{Gi95,GGG03}:
 $$ d [a(\zeta), \zeta, d\zeta ^{\{n-2\}}]=-\frac 1
 {n-2}(\sum_{1\leq j\leq n} \frac {\partial a_j(\zeta)}{d\zeta_j})[\zeta,d\zeta^{\{n-1\}}].$$
 The application of this formula for the column-function
 $a(\zeta)=\frac {u+x}{(\zeta \cdot (u-x))^{n-1} }$ shows the
 closeness of $\kappa$ on $\zeta$ (since $\Delta (u)=\Delta (x)=1$);
  it is closed on $u$ since $[u,du^{\{(n-1\}}]$
 has the maximal degree on $S$.

 The form $\kappa$ is the source of inversion formulas. We
 integrate it on $\Gamma(x)=S\times \gamma(x)$ where $\gamma (x)$ is a cycle
 in
 $\Bbb C^n_\zeta$. Let us now define $\hat f(\zeta, p)$ for any $(\zeta, p)$
 such that the intersection of $S$ by the hyperplane $\zeta
 \cdot z=p$ is empty (we temporarily removed the condition $\Delta
 (\zeta)=0$). If  $(\zeta, p=\zeta\cdot x)$ in the integral satisfy
 these conditions, then the integral makes sense and we can
 interpret the integrand as a differential operator of $\hat f
 (\zeta, p)$.

Let us start from the cycle $\gamma_0(x)$ of real $\zeta=\xi$ such
that $\xi \cdot x=0$. We will pick up the cycle such that it
intersects each ray of 0 in one point (we can take as $\gamma_0(x)
$ the intersection of the sphere $\Delta (\xi)=1$ by the
hyperplane $\xi \cdot x=0$). Of course the integral does not
depend on the specifics of such a choice. Let us restrict $\kappa$ on the
corresponding cycle $\Gamma_0(x)$ and we will interpret it as the
form

$$\kappa_x[f]|_{\Gamma_0}=\frac {f(u)}{(\xi \cdot u-x-i0)^{n-1}}
[u,du^{\{(n-1\}}]\wedge [x+u,\xi,d\xi ^{\{n-2\}}].$$ Here for the
regularization we use the distribution $(\xi\cdot u -i0)^{-n-1}$.
Let us present this form as the sum of 2 forms $\kappa_1+\kappa_2$
corresponding to $x,u$ in $x+u$.

If $f$ is even on $S$ then we have

$$ \int_{ \Gamma_0(x)}
  \kappa_1=\frac {2(2\pi i)^{n-1}}{(n-1)!}f(x).$$

It is just the inversion of the Minkowski-Funk transform
(integration of even functions on $S$ along sections by $\xi
\cdot u=0$). This inversion formula is the specialization of the
inversion of the projective Radon transform \cite {GGG03} (if one were to
present the projective space by pairs of antipodal points) and
equivalent to the inversion of the affine Radon transform.

This integral is equal to zero for odd functions; the
corresponding integral of $\kappa_2$ is equal to zero for even
functions, but will reproduce the odd ones. To see it (using the
rotation invariancy) it is sufficient to consider one point
$x=(1,0.\dots.0)$. Then on $\gamma_0(x)$ we have $\xi_1=0$ and the
integral of $\kappa_2$ for an odd function $f$ coincides with the
integral of $\kappa_1$ for the even function $u_1f(u).$ As the
result, we have for all functions on $S$

$$ \int_{ \Gamma_0(x)}
  \kappa_x [f]=\frac {2(2\pi i)^{n-1}}{(n-1)!}f(x).$$

Now we want to integrate $\kappa$ on the cycle $\Gamma_1(x)$
corresponding to the cycle $\gamma_1(x)$ of horospheres passing
through $x$ (see above cycles on $\partial \Xi _+$). First, we
remark that the form $\kappa$ on $\zeta$ with $\Delta(\zeta)=0$
(corresponding to horospheres) can be transformed in

$$\tilde \kappa_x[f]=\frac {f(u)(\zeta\cdot(u+x))}{(\zeta \cdot x)(\zeta \cdot
(u-x))^{n-1}} [u,du^{\{(n-1\}}]\wedge [x,\zeta,d\zeta
^{\{n-2\}}],\qquad \Delta(\zeta)=0.$$ To see it let us recall that
on the cone $\Xi$ if ${\lambda \cdot
\zeta}\neq 0$  then the form
$$\frac {[\lambda ,\zeta,d\zeta ^{\{n-2\}}]} {\lambda \cdot
\zeta}$$ is independent of $\lambda$. This is simple to see by transforming the determinant
using the condition $\Delta (\zeta)=0$, but the reason for this a
phenomenon is that this form, up to a constant factor, is the
residue of the form $[\zeta,d\zeta ^{\{n-1\}}]/\Delta(\zeta)$ on
$\Xi$ for any $\lambda$. A special case of this independence gives
$$[x+u,\zeta,d\zeta ^{\{n-2\}}]=\frac {\zeta\cdot (x+u)}{\zeta \cdot
x}[x,\zeta,d\zeta ^{\{n-2\}}]$$ and we transformed the form
$\kappa$ in $\tilde \kappa$. This transformation is important,
since we want maximally eliminate $u$ out of the integrand. Now
$u$ participates only trough $p=\zeta \cdot u$.

The direct computation gives that $$ \frac 1 {\zeta \cdot x}\int_S
\frac {f(u)(\zeta\cdot (u+x))}{(\zeta \cdot (u-x))^{n-1}}
[u,du^{\{(n-1\}}]=\frac {(-1)^{n-3}}{(n-2)!}\Cal L_p \hat
f(\zeta,\zeta\cdot x)$$ and since on the cycle $\gamma_1(x)$  we
have $\zeta=x+i\eta, x\cdot \eta=0, \zeta \cdot x=1$, we obtain

$$\int_{\Gamma_1(x)}\kappa_x[f]=i^{n-1}\int_{\gamma_1(x)}\Cal L_p
\hat f(x+i\eta,1)[x,\eta,d\eta^{n-2}].$$

 To finish the proof
of Theorem it is sufficient to prove that the inversion formula
for Radon's cycle $\Gamma_0(x)$ holds also for the horospherical
cycle $\Gamma_1(x)$. Since the form $\kappa$ is closed it is
sufficient to construct a homotopy  of the cycle $\gamma_0(x)$ to
the cycle $\gamma_1(x)$ through such $\zeta$ that the section of
$\Bbb CS $ by the hyperplane $\zeta \cdot z=\zeta \cdot
x-i\varepsilon$ has no real points.

It is very simple to construct such a deformation explicitly.
Using the invariancy it is sufficient to consider only one point
x. Let $x=(1,0\dots,0)$. Then for $\xi \in \gamma_0(x)$ we have
$\xi_1=0$; let $\xi=(0,\tilde \xi)$. Consider cycles
$$\gamma_\delta=\{\zeta(\tilde \xi,\delta)=(i\delta \sqrt{\Delta(\tilde \xi)},\tilde
\xi)\},\qquad 0\leq \delta \leq 1.$$ We have
$\gamma_0=\gamma_0(x); \gamma_1=\gamma_1(x)$ for this $x$. Then
all complex hyperplanes $\zeta(\tilde \xi)\cdot z=p(\tilde
\xi,\delta)=i \delta \sqrt{\Delta(\tilde \xi)}$ intersect the real
sphere $S$ only on the point $x$ and if one were to replace this $p$ on
$p(\tilde \xi,\delta)+i\varepsilon, \varepsilon >0$ then the
hyperplane will have no intersections at all. As a result, we can
define $\hat f(\zeta(\tilde \xi,\delta),p(\tilde
\xi,\delta),0\leq\delta\leq 1$ as the boundary values when
$\varepsilon \rightarrow 0$. So we constructed the desirable
deformation and Theorem is proved.

The horospherical transform $f\rightarrow \hat f$ commutes with
rotations (the group $SO(n)$). The action of $SO(n)$ on the cone
$\Xi$ commutes with the action of $\Bbb C^{\times}$:$\zeta\mapsto
\lambda \zeta$. The domain $\Xi_+$ is not invariant relative to
this action, but the action of the circle $T=\{|\zeta|=1\}$
preserves $\Xi_+$. The decomposition of holomorphic functions on
$\Xi_+$ in Fourier series on $T$ corresponds to the decomposition
on subspaces of homogeneous polynomials and it is the
decomposition on irreducible representations of $SO(n)$:

\proclaim {Proposition} The horospherical transform $f\rightarrow
\hat f$ commutes with the actions of $SO(n)$ and intertwines
subspaces of spherical polynomials on $S$ and homogeneous
polynomials on $\Xi$.\endproclaim

This representation of spherical polynomials through homogeneous
polynomials on the complex cone $\Xi$ goes back to Maxwell's
formulas for spherical polynomials {\cite {B53}, p.251}. This
inverse horospherical transform is in a sense the generating
function for Maxwell's formula. Several components of this
construction admit generalizations. If one were to replace $S$ by
any surface $\{\Phi(x)=1\}$ where $\Phi$ is a polynomial and
$\Phi(u)-\Phi(x)=\sum (u_j-x_j)\phi_j(u,x)$ then the form

$$\kappa_x[f]=\frac {f(u)}{(\zeta \cdot (u-x))^{n-1}}[u,du^{\{n-1\}}]\wedge [\phi(u,x),\zeta,d\zeta ^{\{n-2\}}]$$ is closed, but
it is unclear how to build interesting cycles different from Radon's
cycle. For some examples with non definite quadratic forms cf. in
\cite {Gi00}.

The most natural development of this example is the construction
of horospherical transform on arbitrary compact symmetric spaces
using complex horospheres which is connected with harmonic
analysis on such spaces through Fourier series. This will be the
subject of our future publication.


\Refs \widestnumber \key {GGG03}

\ref \key {B53} \book Higher Transcendental Functions. Bateman
manuscript project \ed Erd\'elyi \vol II \yr 1953  \publ
McGraw-Hill \endref

 \ref \key {GGG03} \by I.Gelfand, S.Gindikin, M.Graev \book
Selected Topics in Integral Geometry \publ Amer.Math.Soc. \yr 2003
\endref

\ref \key {Gi95} \by S.Gindikin \paper Integral geometry on real
quadrics \jour Amer.Math Soc.Transl.(2) \vol 169 \yr 1995 \pages
23-31 \endref

\ref \key {Gi96} \by S.Gindikin \paper Real integral geometry and
complex analysis \inbook Integral Geometry, Radon Transform and
Complex Analysis. Lect.Notes in Math. \vol 1684 \yr 1998 \pages
70-98 \publ Springer \endref

\ref \key {Gi00} \by S.Gindikin \paper Integral Geometry on
$SL(2;\Bbb R)$ \jour Math.Res.Letters \yr 2000 \pages 417-432
\endref
\endRefs
\enddocument